\documentclass{amsart}

\usepackage{rlepsf, amssymb, amsmath, amsthm}

\newtheorem{thm}{Theorem}[section]
\newtheorem{lem}[thm]{Lemma}

\newtheorem{prop}[thm]{Proposition}
\newtheorem{rem}[thm]{Remark}
\newtheorem{example}[thm]{Examples}

\newcommand{\bbr}{\begin{rem}\em} 
\newcommand{\eer}{\end{rem}}

\newcommand{\bex}{\begin{example}\em} 
\newcommand{\eex}{\end{example}}

\def\bqtm#1{\begin{proof}[{\em \bf Theorem~\ref{#1}}] \em}
\def\eqtm{\end{proof}}

\def\dfn#1{{\em #1}}

\def\co{\colon\thinspace}

\def\Z{\text{$\mathbb{Z}$}}

\def\R{\text{$\mathbb{R}$}}
\def\C{\text{$\mathbb{C}$}}

\begin{document}
\title{Planar open book decompositions and contact structures}

\author{John B.\ Etnyre}
\address{Department of Mathematics,
University of Pennsylvania,
209 South 33rd St.,
Philadelphia, PA 19104-6395}
\email{etnyre@math.upenn.edu}
\urladdr{http://math.upenn.edu/\char126 etnyre}

\begin{abstract}
In this note we observe that while all overtwisted contact structures on compact 3--manifolds are
supported by planar open book decompositions, not all contact structures are. This has relevance to
the Weinstein conjecture \cite{AbbasCieliebakHofer} and invariants of contact structures.
\end{abstract}
\maketitle

\section{Introduction}
In \cite{Giroux??}, Emmanuel Giroux showed that all contact structures on compact 3--manifolds come
from open book decompositions via the Thurston-Winkelnkemper \cite{TW} construction. That is, given any contact
structure $\xi$ on $M$ there is an open book decomposition of $M$ such that $\xi$ is transverse to
the binding of the open book and can be isotoped arbitrarily close to the pages. This
fundamental breakthrough has provided the basis for a much greater understanding of contact
structures, see \cite{Eliashberg04, Etnyre04, EtnyreHonda02a,OzsvathSzabo:Contact}, and 3--manifold
topology, see \cite{KronheimerMrowka, OzsvathSzabo}.

A first obvious question concerning the connection between open book decompositions and contact
structures is what is the minimal genus of a page of an open book decomposition that supports a
given contact structure. It is interesting to observe that {\em any} 3--manifold has an open book
decomposition with planar pages (that is, the pages are $S^2$ with a finite number of disjoint open disks
removed) \cite{Rolfsen}. Even given this it seems unlikely that every contact structure can be
supported by a planar open book. However, for overtwisted contact structures we have:
\bqtm{otandplanar}
Any overtwisted contact structure on a closed 3--manifold is supported by a planar open book
decomposition.
\eqtm
Thus there is no homotopy theoretic obstruction to a contact structure admitting a planar open book, but there
are nonetheless obstructions. To state them we first recall that given a 4--manifold with boundary
$X$ we can consider its intersection form $q_X$ on $H_2(X,\Z).$ We denote by $b_2^\pm$ the maximal
dimension of a subspace of $H_2(X,\Z)$ on which $q_X$ is $\pm$ definite and we denote $b_2^0(X)$ the
dimension of the subspace on which $q_X$ is degenerate. 
\bqtm{planarobstruction} 
If $X$ is a symplectic filling of a contact 3--manifold $(M,\xi)$ supported by a planar open book
decomposition then $b_2^+(X)=b_2^0(X)=0$ and the boundary of $X$ is connected. Moreover, if $M$ is
an integral homology sphere then the intersection form on $X$ is diagonalizable. 
\eqtm
We now give several simple applications of this theorem. 

\bex 
Let $(M,\xi)$ be a contact manifold obtained by performing Legendrian surgery on a Legendrian link
in $S^3$ with one component having Thurston--Bennequin invariant larger than 0, then $(M,\xi)$ is
not supported by an planar open book decomposition.
\eex

\bex 
In \cite{McDuff91}, McDuff gave the first examples of a symplectic manifold $(X,\omega)$ with two
convex boundary components. For example, take the unit disk bundle in the cotangent bundle of a
surface of high genus, perturb the symplectic form so that the 0-section is symplectic, and let
$(X,\omega)$ be the complement of a symplectic neighborhood of this perturbed section. It is not
hard to check that both boundary components of $X$ are convex. Let $(M,\xi)$ be either boundary
component $(X,\omega)$ then $(M,\xi)$ cannot be supported by a planar open book.
\eex

\bex 
Let $(X,\omega)$ be the symplectic manifold obtained by plumbing together symplectic disk bundles
over $S^2$ with Euler number $-2$ according to the $-E_8$ graph. The boundary of $X$ is convex and
topologically the Poincar\'e homology sphere $P.$ Let $\xi$ be the tight contact structure induced
on $P$ by $(X,\omega).$ According to Theorem~\ref{planarobstruction}, $(P,\xi)$ cannot be supported
by a planar open book decomposition since the intersection from of $X$ is non-diagonalizable.
\eex

It would be very interesting to find obstructions to a contact structure being supported by a planar
open book that were more intrinsically 3--dimensional. In particular, an obstruction that could be
used to show some non-fillable contact structures are not supported by planar open books. It would
also be interesting to find bounds on the minimal genus of a supporting open book for a contact
structure. It is still conceivable that all contact structures are supported by open book
decompositions with genus one pages.

Given a contact form $\alpha$ for the contact structure $\xi$ on $M$ the \dfn{Reeb vector} field
$v_\alpha$ is the unique vector field satisfying $\alpha(v_\alpha)=1$ and
$\iota_{v_\alpha}d\alpha=0.$ Recall the well known (extension of the) Weinstein conjecture
\cite{Weinstein79} asserts any Reeb vector field $v_\alpha$ for $\xi$ has a closed periodic
orbit. This conjecture was proven of any contact structure on $S^3$ or on any reducible 3--manifold
by Hofer in \cite{Hofer93}. That paper also establishes the conjecture for all contact structure
that are virtually overtwisted (that is have a finite cover that is overtwisted). Abbas, Cieliebak
and Hofer have a program for proving the Weinstein conjecture based on open book
decompositions. Currently \cite{AbbasCieliebakHofer} they can prove the Weinstein conjecture is true
for any contact structure supported by a planar open book decomposition. While this greatly enlarges
the class of contact structures for which the Weinstein conjecture is known the above examples and
Theorem~\ref{planarobstruction} show there are contact structures for which the conjecture is still not
yet known.

In Section~\ref{background} we recall the necessary background about open book decompositions,
contact surgeries and symplectic surface bundles. In the following section we establish
Theorem~\ref{otandplanar} concerning open books for overtwisted contact structures. In
Section~\ref{sec:obstop} we derive obstructions for a contact structure to be supported by a planar
open book. 

Acknowledgments: I am grateful to Kai Cieliebak and Helmut Hofer for their inspiring questions
about open books and contact structures and their encouragement with this project. I thank Ko
Honda for several enlightening conversations and helpful comments on the paper. I am also grateful
to Burak Ozbagci who provided many helpful comments on the first draft of this paper and to Andy Wand
for pointing out an error in the original version of Figure 3. In addition, I 
gratefully acknowledge partial support by an NSF Career Grant (DMS-0239600) and FRG-0244663.

\section{Open book decompositions and contact structures}\label{background}
In this section we recall various facts about open book decompositions, contact surgeries and
symplectic surface bundles.

\subsection{Open book decompositions and contact strucutres}
Recall an \dfn{open book decomposition} of a 3--manifold $M$ is a triple $(B, \Sigma,\phi)$ where $B$ is a link in 
$M$ such that $M\setminus B$ fibers over the circle with fiber $\Sigma$ and monodromy $\phi$ so that $\phi$ is the identity
near the boundary and each fiber of the fibration is a Seifert surface for $B.$ By saying $\phi$ is
the monodromy of the fibration we mean that $M\setminus B = \Sigma\times [0,1]/\sim$ where
and $(1,x)\sim (0,\phi(x)).$ The fibers in the fibration are called \dfn{pages} of the open
book and $B$ is called the \dfn{binding}. Note, given a diffeomorphism of a surface $\phi$ that is
fixed near the boundary we may form it's mapping torus and glue in solid tori to get a closed three
manifold having an open book decomposition $\phi.$ So sometimes we will designate an open book
decomposition simply by $\phi$ without reference to the binding $B.$ (The main difference here is
whether we are thinking of the open book as inside a pre-existing 3--manifold or whether we are
defining the 3--manifold by the open book.)

Two open books for a manifold $M$ are said to be \dfn{equivalent} if there is an ambient isotopy of $M$
taking the binding and pages of one to the binding and pages of the other. Given an open book
$(B, \Sigma, \phi)$ for $M,$ let $\Sigma'$ be $\Sigma$ with a 1--handle attached. Let
$\alpha$ be an simple closed curve in $\Sigma'$ that intersects the co-core of the attached
1--handle exactly once. Set $\phi'=\phi\circ D^\pm_\alpha,$ where $D^\pm_\alpha$ is a right/left
handed Dehn twist along $\alpha.$ The mapping torus of $\phi'$ has torus boundary components
each with a canonical product structure $S^1\times S^1$ where the first $S^1$ bounds the fibers in
the fibration of the mapping torus. Let $M'$ be the mapping torus with $S^1\times D^2$ glued to each
boundary component so that the product structure (and ordering of the $S^1$ factors) is
preserved. Let $B'$ be the cores of the added tori. One may check that $M'$ is diffeomorphic to $M$
and thus $(B', \Sigma', \phi')$ is another open book for $M.$ If $D^+_\alpha$ was used $(B', \Sigma', \phi')$ is called
the \dfn{positive stabilization} of $(B,\phi)$ otherwise it is called the \dfn{negative stabilization}.

A contact structure $\xi$ on $M$ is \dfn{compatible with}, or \dfn{supported by} an open book $(B, \Sigma, \phi)$ of $M$ if 
$B$ is transverse to $\xi$ and on the complement of $B$ the contact planes $\xi$ can be isotoped to
be arbitrarily close to the pages of the open book (while keeping $B$ transverse). Thurston and
Winkelnkemper \cite{TW} showed that any open book supports a contact structure. It is not too hard
to see that if two contact structures are supported by an fixed open book then they are isotopic as
contact structures. One may also check that if $(B, \Sigma, \phi)$ is compatible with a contact structure
then so is any open book obtained from $(B,\Sigma,\phi)$ by {\em positive} stabilization. 
Giroux made the following fundamental observation (for a discussion of the proof
see \cite{EtnyreOBN, Goodman}).
\begin{thm}[Giroux \cite{Giroux??}]
Every contact structure is supported by an open book. Moreover, if two open books support the same
contact structure then they each may be positively stabilized some number of times so that the open
books are equivalent.
\end{thm}

Suppose $(B_1, \Sigma_1, \phi_1)$ and $(B_2,\Sigma_2, \phi_2)$ are open book decompositions supporting the contact manifolds
$(M_1,\xi_1)$ and $(M_2,\xi_2),$ respectively. Let $\alpha_i$ a properly embedded arc in $\Sigma_i.$ 
The Murasugi sum of $(B_1, \Sigma_1, \phi_1)$ and $(B_2, \Sigma_2, \phi_2)$ is obtained as follows: let $\Sigma_1 * \Sigma_2$
be the surface $\Sigma_1\cup_{R_1=R_2} \Sigma_2$ where $R_i$ is a rectangular neighborhood of $\alpha_i$ in $\Sigma_i$ and
$R_1$ is identified to $R_2$ in such a way that 
$\partial R_1\cap \partial \Sigma_1=\overline{\partial R_2\setminus \partial \Sigma_2}.$
Each of the $\phi_i$'s can be extended to $\Sigma_1 * \Sigma_2.$ The Murasugi sum in now defined to be the open book obtained
from $\phi_1\circ \phi_2.$ 
\begin{lem}\label{murasugi}
If $(B_1, \Sigma_1, \phi_1)$ and $(B_2,\Sigma_2, \phi_2)$ are open book decompositions supporting the contact manifolds
$(M_1,\xi_1)$ and $(M_2,\xi_2),$ respectively, then the Murasugi sum of $(B_1, \Sigma_1, \phi_1)$ and $(B_2,\Sigma_2, \phi_2)$
is an open book decomposition supporting the contact manifold $(M_1,\xi_1)\#(M_2,\xi_2).$
\end{lem}

\subsection{Contact Surgeries}
Let $(M,\xi)$ be a contact manifold and $L\subset M$ a closed Legendrian curve.  Let 
$N(L)$ be a {\em standard tubular neighborhood} of the Legendrian curve $L.$ This means the
neighborhood has convex boundary and two parallel dividing curves (see \cite{EtnyreHonda01b}).  Choose a framing for $L$  
so that the meridian has slope $0$ and the dividing curves have slope $\infty$.  With respect to 
this choice of framing, a \dfn{$\pm 1$ contact
surgery} is a $\pm 1$ Dehn surgery, where a copy of $N(L)$ is glued to $M\setminus N(L)$ so that the new 
meridian has slope $\pm 1$.  Even though the boundary characteristic foliations may not 
exactly match up a priori,  we may use Giroux's Flexibility Theorem 
\cite{Giroux91, Honda00} and the fact that they have the same dividing set to make the 
characteristic foliations agree. This gives us a new manifold $(M',\xi').$ It is common to call $-1$
contact surgery \dfn{Legendrian surgery}. For a detailed discussion of contact surgery see
\cite{DingGeigesStipsicz04}. The following is a well known theorem, see for example \cite{EtnyreHonda02a}.
\begin{thm}\label{LegSurgonPage}
Suppose the $L$ is a Legendrian knot in the contact manifold $(M,\xi), \xi$ is supported by the open
book $(B,\phi)$ and $L$ is contained in a page of the open book. The contact manifold
obtained from $(M,\xi)$ by $\pm 1$ contact surgery on $L$ is equivalent to the one compatible with
the open book with monodromy $\phi\circ D^\mp_\alpha.$
\end{thm}

\subsection{Symplectic surface bundles}
In this paper a \dfn{symplectic fibration over $S^1$} will mean a 3--manifold $M$ that fibers over the
circle together with a closed 2-form $\omega$ which is positive on each fiber. The kernel of
$\omega$ defines a line field $l$ that is transverse to the fibers of the fibration. An orientation
on $M$ and on the fibers induces an orientation on $l$ and thus we can fix a fiber $\Sigma_0$ of the
fibration and use $l$ to define a return map $H_{(M,\omega)}\co \Sigma_0\to \Sigma_0$ called the
\dfn{holonomy} of the symplectic fibration. If we normalize $\omega$ so that it integrates to 1 on
each fiber of the fibration then the holonomy determines $(M,\omega)$ up to fiber preserving
diffeomorphism. It is also important to notice that the holonomy determines a symplectic
neighborhood of $(M,\omega)$ or more precisely we have the following lemma.
\begin{lem}\label{holdet}
Suppose $(X_i,\omega_i)$ is a symplectic 4--manifold and $f_i:M\to X_i$ is an embedding of a fibered
3--manifold into $X_i$ for $i=0,1.$ Suppose $f_i^*\omega_i$ define symplectic fibrations on $M$ with
the same holonomy. Then $f_0(M)$ and $f_1(M)$ have symplectomorphic neighborhoods.
\end{lem}
We will also find the following lemma useful.
\begin{lem}\label{confill}
Suppose the holonomy of the symplectic fibration $(M,\omega)$ is a Hamiltonian diffeomorphism.
Then there is a symplectic form $\Omega$ on $X=\Sigma\times D^2,$ such that $\partial
(X,\Omega)=(M,\omega),$ where $\Sigma$ is the fiber of the fibration. 
\end{lem}

We now have the fundamental result of Eliashberg \cite{Eliashberg04}.
\begin{thm}\label{collar}
Suppose $(M,\xi)$ is a contact 3--manifold and $\omega$ is a closed 2--form on $M$ that is positive
on $\xi.$ Furthermore suppose $(B,\phi)$ is an open book supporting $\xi.$ Use the pages of the open
book to define a framing on the components of the binding $B.$ Let $X$ be the 4--manifold obtained
from $M\times [0,1]$ by attaching 2-handles to $M\times\{1\}$ along $B$ with framing 0. The manifold
$X$ is oriented so that $\partial X = (-M)\cup M',$ where $M'$ is the fibered three manifold
obtained by 0--surgery on $B$ in $M.$ Then $X$ admits a symplectic form $\Omega$ such
that $\Omega|_M=\omega$ and $\Omega|_{M'}$ defines a symplectic fibration on $M'.$ 
\end{thm}

\section{Overtwisted contact structures}
In order to prove all overtwisted contact structures are supported by open books we need a few
preliminary observations.
\begin{lem}\label{genhom}
Suppose $(B,\phi)$ is an open book for the 3--manifold $M.$ Let $\alpha_1,\ldots, \alpha_k$ be
curves on a page of the open book that generate the homology of the page. Then the homology of $M$
is generated by $\alpha_1,\ldots,\alpha_k.$\qed
\end{lem}
The proof of this lemma is an easy exercise that is left to the reader. Recall that given an
oriented Legendrian knot $L$ then there is a positive $L_+$ and negative $L_-$ transverse push off,
\cite{Eliashberg93, EtnyreHonda01b}. Moreover, there is a positive $S_+(L)$ and negative $S_-(L)$ 
stabilization of $L,$ defined by ``adding down or up zig zags''. (This is completely accurate for
Legendrian knots in the standard contact structure on $\R^3$ represented by their front
projection. For Legendrian knots in a more general contact structure one must take a little more
care in the definition \cite{EtnyreHonda01b}.)   
\begin{prop}[Ding, Geiges and Stipsicz \cite{DingGeigesStipsicz}]\label{prop:ots}
Let $L$ be a Legendrian knot in a contact 3--manifold $(M,\xi).$ Let $L'$ be a parallel copy of $L$
that has been stabilized positively (respectively negatively) twice. Let $T$ be $L_+$ (respectively
$L_-$). The contact structure obtained from $\xi$ by performing a Lutz twist along $T$ is
isotopic to the one obtained from $\xi$ by performing $+1$-contact surgery along both $L$ and $L'.$
\end{prop}
The proposition was first hinted at in Proposition~4.3 of \cite{EtnyreHonda02a}. It was explicitly
stated and given a nice proof by Ding, Geiges and Stipsicz \cite{DingGeigesStipsicz}. 
Here we sketch a somewhat different proof.
\begin{proof}[Sketch of Proof]
Note that both surgeries take place in a solid torus. So we will restrict our attention to a
neighborhood of $L$. Let $N$ be a standard neighborhood of a Legendrian knot $L.$ That is $\partial
N$ is convex with two dividing curves of slope $\infty.$ Given that $\xi$ restricted to $N$ is
tight, this uniquely determines a contact structure on $N$ (see \cite{Giroux00, Honda00}). We can
write $N$ as $N_1\cup N'$ where $N_1$ is a solid torus with convex boundary having two dividing
curves of slope $-\frac12$ and $N'=\overline{N\setminus N_1}.$ There are actually three ways of
splitting $N$ like this. They come from stabilizing the core Legendrian knot in $L$ twice
positively, twice negatively or once positively and once negatively. See \cite{EtnyreHonda01b}. Now
split $N'$ into two pieces $N_2\cup N_3$ where $N_3$ is a vertically invariant neighborhood of
$\partial N$ that is contained in $N'$ and $N_3$ is the complement of this neighborhood. One may
easily check that a Legendrian divide on $N_3\cap N_2$ is Legendrian isotopic to $L.$ Also $N_1$ is
a standard neighborhood of a Legendrian knot obtained from $L$ by stabilizing twice. That is, with
the appropriate choice of $N_1,$ $N_1$ is a standard neighborhood of $L'.$ For this and the
following facts about contact structures on toric annuli and Legendrian knots see
\cite{EtnyreHonda01b}.

We can identify slopes $s\in\R\cup\{\infty\}$ of linear foliations on $T^2$ with their respective
angles, $[\theta_s]\in \R/\pi \Z$. In order to distinguish the different amounts of twisting, we
will choose a lift $\theta_s\in \R$ instead. There exists an exhaustion of $N$ by concentric $T^2$
with linear foliations where the angles of the foliation on the tori monotonically increase over the
interval $(\frac{\pi}{2}, \pi)$ as the $T^2$ move towards the core of $N.$

We can think of $L$ as sitting in $N_3$ as a Legendrian divide on a convex torus in $N_3.$ One may
check that performing $+1$-surgery on $L$ gives a tight minimally twisting contact structure on
$N_3.$ Moreover, the dividing curves on both boundary components still have infinite slope. (The
product structure on the $T^2\times I$ has changed though.) The contact structure, after surgery
on $L,$ on $N=N_1\cup N_2\cup N_3$ is tight and we can find concentric $T^2$ where the
characteristic foliation are linear and run over the interval $(\frac{\pi}{8}, \pi).$ Performing
$+1$-surgery on $L'$ yields a contact structure on $N$ with twisting from $(-\frac{\pi}{2},\pi).$
Thus the result of performing the $+1$-contact surgeries on $L$ and $L'$ is equivalent to
replacing the contact structure on $N$ that twists over $(\frac{\pi}{2}, \pi)$ by one that twists
over $(-\frac{\pi}{2},\pi).$ But this is exactly the change that happens when one does a Lutz twist
on a transverse push off of the core of $N.$
\end{proof}

\begin{figure}[ht]
  \relabelbox \small {\epsfxsize=3.5in\centerline{\epsfbox{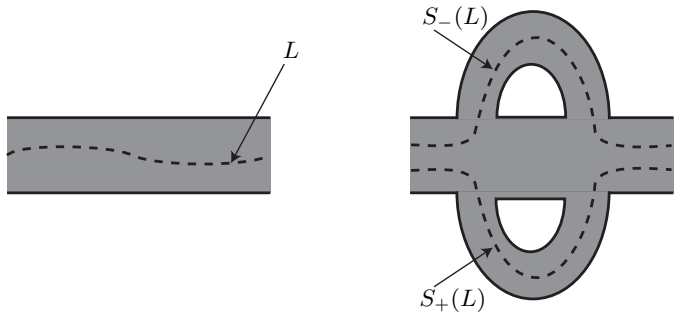}}} 
  \relabel{l}{$L$}
  \relabel {p}{$S_+(L)$}
  \relabel {n}{$S_-(L)$}
  \endrelabelbox
        \caption{A neighborhood of a piece of $L$ in $\Sigma,$ left.($L$ is oriented so it points towards the left.) 
          The twice stabilized open book, right.}
        \label{stab}
\end{figure}

In our discussion below it will be useful to see how to stabilize a Legendrian knot on a page of an open
book so that the stabilized knot is also on a page of an open book. 
\begin{lem}\label{lem:stab}
Let $(B,\Sigma,\phi)$ be an open book decomposition supporting the contact structure $\xi$ on $M.$ Suppose
$L$ is a Legendrian knot in $M$ that lies on a page of the open book. If we positively stabilize 
$(B,\Sigma,\phi)$ twice as shown in Figure~\ref{stab} then we may isotop the page of the open book so that
$S_+(L)$ and $S_-(L)$ appear on the page as seen in Figure~\ref{stab}.
\end{lem}
This Lemma is relatively easy to prove, see \cite{EtnyreOBN}.

We first prove our main theorem for the special case of $M=S^3.$
\begin{lem}
Any overtwisted contact structure on $S^3$ is supported by a planar open book decomposition.
\end{lem}
\begin{proof}
We simply need to exhibit a planar open book for $S^3$ realizing each homotopy class of plane field. Using
the notation of \cite{DingGeigesStipsicz04} ({\em cf} \cite{Gompf98}) the homotopy classes of plane fields correspond to half-integers.
That is given a plane field $\xi$ on $S^3$ let $(X,J)$ be an almost complex 4-manifold with $\partial X=S^3$ and
$\xi$ the set of complex tangencies to the boundary. Then we associate the invariant
\[d_3=\frac14(c_1(X,J)^2 - 3\sigma(X) -2 \chi(X)).\]
One may compute (see \cite{DingGeigesStipsicz04}) 
that the indicated contact surgeries in Figure~\ref{seg} represent two overtwisted contact
structures on $S^3$ with $d_3$ equal to $\frac12$ and $-\frac32,$ respectively. From Lemma~\ref{lem:stab} 
\begin{figure}[ht]
  \relabelbox \small {\epsfxsize=3.5in\centerline{\epsfbox{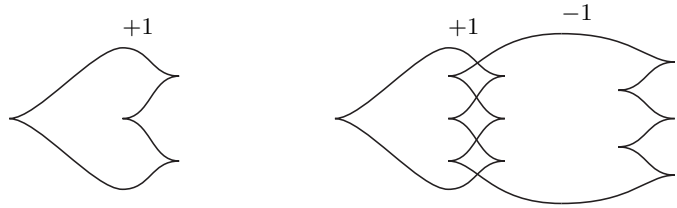}}} 
  \relabel{1}{$+1$}
  \relabel{2}{$-1$}
  \relabel{3}{$+1$}
  \endrelabelbox
        \caption{Surgery diagrams for overtwisted contact structures with $d_3=\frac12,$ left, and $d_3=-\frac32,$ right.}
        \label{seg}
\end{figure}
one may easily see
that the two contact structures in Figure~\ref{seg} are supported by the open books shown in Figure~\ref{ob}. 
\begin{figure}[ht]
  \relabelbox \small {\epsfxsize=3.5in\centerline{\epsfbox{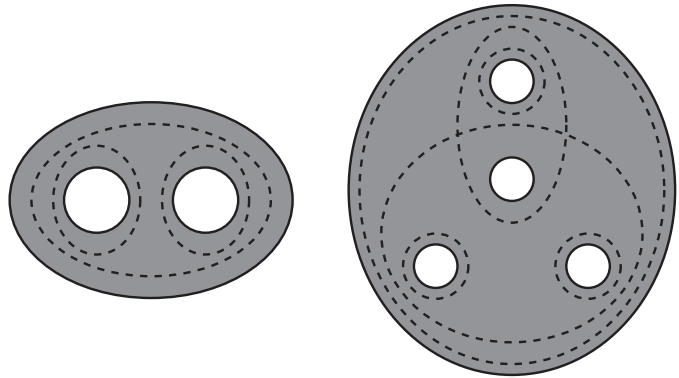}}} 
  \endrelabelbox
        \caption{The fibers of open book decompositions for the contact structures in Figure~\ref{seg}.
        The monodromies come from Dehn twists about the dotted curves in the picture. All Dehn twists are
      right handed except for the Dehn twist about the outer most curve in both pictures which is left handed.}
        \label{ob}
\end{figure}
Note these
are both planar open books. The invariant $d_3$ is additive under contact connected sum, so by taking the
connected sum of these examples we can realize any overtwisted contact structure on $S^3.$ These connected sums
can be obtained by Murasugi summing the corresponding open books (see Lemma~\ref{murasugi}). Since Murasugi 
sums can be performed in a way so that the genus of the fibers adds, we see that we have planar open books supporting
all overtwisted contact structures on $S^3.$
\end{proof}

We are now ready for the main theorem of this section.
\begin{thm}\label{otandplanar}
Let $(M,\xi)$ be an overtwisted contact 3--manifold. Then there is an open book $(B,\phi)$ for $M$
that supports $\xi$ and has planar pages.
\end{thm}
\begin{proof}
Two plane field are homotopic if and only if their ($\frac12$) Euler classes $d_2$ and their three dimensional 
invariants $d_3$ are the same (see \cite{GeigesIntro} for a description of this using the notation we use here, 
also see \cite{Gompf98}). 
We show that all such invariants can be realized
on a manifold by overtwisted contact structures supported by planar open books. Thus Eliashberg \cite{Eliashberg89} 
tells us that all overtwisted contact structures will then be realized by planar open books.

It is well known, see \cite{Rolfsen}, that any three manifold $M$ has an open
book with planar leaves $(B,\Sigma,\phi).$ The associate contact structure $\xi$ may or may not be overtwisted, but all overtwisted
contact structures can be obtained from it.
We begin by showing that any possible $d_2$ can be realized by an overtwisted contact structure supported by a planar
open book. To this end let $\alpha_1,\ldots,\alpha_k$ be the 
simple closed curves on $\Sigma$
that generate the homology of the page. We know that the $\alpha_i$ generate the homology of $M$ by
Lemma~\ref{genhom}.
If we orient $\alpha_i,$ make it transverse (respecting this orientation) and perform a Lutz twist on $\alpha_i$ then one 
may compute that the difference 
between $d_2$ for the two contact structures is Poincar\'e dual to the homology class given by $\alpha_i.$
Thus by performing Lutz twists along the $\alpha_i$'s with various orientation (and transverse realizations)
we can realize any possible $d_2$ by an overtwisted contact structure. 
Using Proposition~\ref{prop:ots}
we may obtain all invariants $d_2$ by overtwisted contact structures by performing various contact surgeries on Legendrian
realizations of the $\alpha_i$'s and their stabilizations. We can Legendrian realize any $\alpha_i$ on 
the page of our open book. To realize all the possible stabilization of $\alpha_i$ that we will need
we must positively stabilize the open book four times for each $\alpha_i.$ (Two stabilizations
of the open book to realize $S_+^2(\alpha_i)$ and two to realize $S_-^2(\alpha_i).$) We can do this keeping
the open book planar. Now all the Legendrian knots on which we need to do contact surgeries can be 
Legendrian realized on the pages of the open book. Performing a $\pm1$ contact surgery on one of these
Legendrian knots will yield a new open book whose monodromy is the old monodromy composed with 
a $\mp$ Dehn twist along the knot in the page. Moreover, this open book supports the contact structure
obtained from $\xi$ by the contact surgeries which were performed. Thus we have realized all possible invariants $d_2$
by overtwisted contact structures supported by planar open books.

Now recall that the invariant $d_3$ is additive under contact connected sum. Thus to realized any pair $(d_2,d_3)$
by an overtwisted contact structure supported by a planar open book we first take $(M,\xi')$ realizing the
appropriate $d_2$ invariant. Then there is a contact structure $\xi''$ on $S^3$ such that $\xi'\# \xi''$ on $M$ 
realizes $(d_2, d_3).$ Since $\xi''$ and $\xi'$ are both supported by planar open books, we can perform the 
Murasugi sum in such a way that the resulting open book is also planar. Thus we have a planar open book 
supporting the contact structure   $\xi'\# \xi''.$
\end{proof}

\section{The genus of open books}\label{sec:obstop}

\begin{thm}\label{planarobstruction}
If $X$ is a symplectic filling of a contact 3--manifold $(M,\xi)$ supported by a planar open book
decomposition then $b_2^+(X)=b_2^0(X)=0$ and the boundary of $X$ is connected. Moreover, if $M$ is an
integral homology sphere then the intersection form on $X$ is diagonalizable.
\end{thm}

\begin{proof}
Let $(X,\omega)$ be a weak filling of $(M,\xi)$ and assume $\xi$ is supported by the planar open
book $(B,\Sigma,\phi).$ For now assume the boundary of $X$ is connected. Theorem~\ref{collar} and
Lemma~\ref{holdet} say we may attach 2--handles to $X$ along $B\subset M,$ with 0 framing measured
with respect to the fibers of the open book, so that the new 4--manifold $X'$ has a symplectic
structure $\omega'$ such that $\omega'|_X=\omega$ and the new boundary component of $X'$ is
$M'=S^2\times S^1$ and $S^2\times\{p\}$ is symplectic for each $p\in S^1.$ That is $\omega'$ gives
$S^2\times S^1$ the structure of a symplectic fibration.

Let $C$ be $S^2\times D^2$ and $\omega_c$ the symplectic form on $C$ coming from Lemma~\ref{confill}
such that $\partial (C,\omega_c)=-(M',\omega'|_{M'}).$ (Note we can apply the lemma since any
symplectomorphism of $S^2$ is Hamiltonian.) Thus using Lemma~\ref{holdet} we can construct the
closed symplectic 4--manifold $(W,\Omega)=(M',\omega')\cup (C,\omega_c).$ 
Clearly $S=S^2\times\{p\}$ is a symplectic sphere in $W$ with square 0. Thus a theorem of
McDuff \cite{McDuff90} implies that $(W,\Omega)$ is the blow up of a ruled surface. The cocore of a
2-handle $H$ attached to $X$ in forming $X'$ can be glued to the $\{q\}\times D^2$ in $C$ to form a
topological sphere $S'$ with square $k$ and intersecting $S$ geometrically one time. If $B_1$ is the
component of the binding $B$ that $H$ is attached to then $k$ depends on the Dehn twists in the 
monodromy that are parallel to $B_1.$ Thus by adding positive Dehn twists to the monodromy if
necessary we can assume that $k$ is even. Now there is a sphere $S''$ in the homology class $S'-\frac{k}{2}S$
that has self intersection 0 and intersects $S$ exactly one time. Thus a
neighborhood $N$ of $S\cup S''$ is a punctured $S^2\times S^2$ in $W.$  McDuff's result now
implies that $W= S^2\times S^2\#_n \overline{\C P}^2.$
Moreover $X\subset W$ is disjoint from $N$ so we may embed $X$ in $\#_n \overline{\C P}^2$ and
since this manifold has a negative definite, diagonalizable intersection form we
see that $b_2^+(X)=b_2^0(X)=0$ and if $M$ is a homology sphere then the intersection from on $X$ is
diagonalizable.

Now if the boundary of $X$ was not connected then we can find symplectic caps for the other boundary
components \cite{Eliashberg04, Etnyre04}. That is, for each extra boundary component $M_i$ of $X$
we can construct a symplectic manifold $(C_i, \omega_i)$ such that $\partial C_i=-M_i$ and
$(C_i,\omega_i)$ may be glued to $(X,\omega)$ along $M_i$ to form a new closed symplectic manifold
with one less boundary components. Moreover, examining the construction in \cite{Etnyre04} one can
easily arrange that $b_2^+(C_i)>0.$ Thus after capping off the extra boundary components of
$(X,\omega)$ we have a filling of $(M,\xi)$ with one boundary component and $b_2^+>0$ contradicting
the first part of our proof. 
\end{proof}

\bibliographystyle{alpha}

\end{document}